# Some Symmetric Pyramids with trivial Dehn invariant

BY G. DUVAL

guillaume.duval@insa-rouen.fr

**Abstract**

For some symmetric pyramids of $\mathbb{R}^3$, we find Galois obstruction for their Dehn invariant to be zero, i.e. for the pyramids to be scissor equivalent to a cube. These conditions are that some associated Kummer extensions of number fields must be abelian. This work can be viewed as a complement to [5]. Indeed, in [5], an impressive classification of all rational tetrahedron is given, while we concentrate our attention on much more symmetric pyramids and prove that the only ones which are scissor equivalent to a cube are the rational ones.

## 1 Symmetric pyramids scissor equivalent to cubes

In [3], Conway and Jones classified all minimal vanishing sums of roots of unity with at most nine terms. One of their geometric motivations was to find possible rational euclidian tetrahedrons of $\mathbb{R}^3$, whose six dihedral angles are rational multiples of $\pi$. Such a performance was achieved only recently in [5]. In their paper, Kedlaya, Kolpakov, Poonen and Rubinstein, manage to classify all rational tetrahedrons of $\mathbb{R}^3$. To this aim, they improved Conway and Jones results by classifying vanishing sums of roots of unity with at mots twelve terms. They where also led to use other involved technics among other the use of the Regge symmetries on the set of all tetrahedrons. According to [1] and [2], a rational tetrahedron is scissor equivalent to a cube since its Dehn invariant is zero. But conversely, there exist euclidian solids which are not rational but which are scissor equivalent to cubes. Since the problem of classifying all solid even the convex ones with trivial Dehn invariant is too difficult, in this article, we will focus our attention on certain symmetric pyramids.

In the left part of Figure 1, a cube is the union of six isometric "egyptian pyramids". Therefore, $\mathrm{Dehn}(\mathrm{Cube}) = 0 = 6\,\mathrm{Dehn}(\mathrm{Pyramid})$ which implies that the Dehn invariant of the pyramid is zero. In other words, the pyramid is scissor equivalent to a cube.

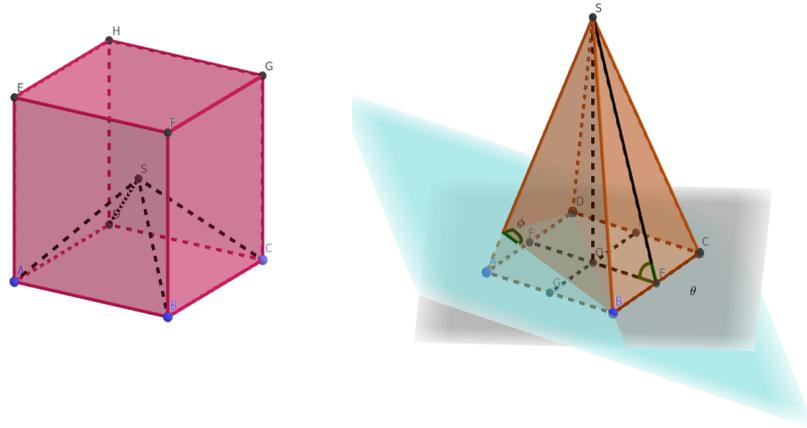

**Figure 1.**

Let $P_4(h)$ be the pyramid of the right part of the figure, where $h = OS$ denote its height and whose square basis is made of the four quartic roots of the unity : $\mu_4 := \{\pm 1; \pm i\}$. For $h = 1/\sqrt{2}$, $P_4(1/\sqrt{2})$ is one of the six pyramid inscribed into the cube. What's happen if we change the height of this pyramid as in the right part of the figure ? More generally, for any real positive number $h > 0$, for any integer $n \geqslant 3$, let's denote by $P_n(h)$ the symmetric pyramid whose horizontal basis is made of the $n$ roots of the unity $\mu_n$ and whose vertical height is $h = OS$. Our main result will be :

**Theorem 1.** *For any $h > 0$, the only right square pyramid $P_4(h)$ with trivial Dehn invariant is $P_4(1/\sqrt{2})$. A a consequence, this pyramid is the only one which is rational, in the sens that all its dihedral angles are rational multiple of $\pi$.*

*When the square basis is changed by equilateral triangles or by regular hexagons, the corresponding pyramids $P_3(h)$ and $P_6(h)$ have non trivial Dehn invariants. In other words they never are scissor equivalent to a cube.*

For the moment let's derive some consequences of the present result, in the square basis case. Similar consequences for $P_3(h)$ and $P_6(h)$ are left to the reader.

If a polyhedron $\mathcal{P}$ can be cut into, two symmetric sub polyhedron : $\mathcal{P} = \mathfrak{P} \bigcup \mathfrak{P}'$ by some plane of symmetry, then $\mathrm{Dehn}(\mathcal{P}) = 2\,\mathrm{Dehn}(\mathfrak{P})$. Since the multiplication by any non zero integer is an injectiv mapping of the Dehn group : $\mathbb{R} \otimes_{\mathbb{Z}} \mathbb{R}/\pi\mathbb{Z}$, we have that $\mathcal{P}$ and $\mathfrak{P}$ both are scissor equivalent to a cube in the same time or not. The next figure shows the six successive symmetric cutting of $P_4(h)$ by vertical planes containing the vertical line $(O\,S)$ and symmetrically cutting the square basis. The resulting operations give six polyhedron $T_1(h) := P_4(h), T_2(h), ..., T_6(h)$. Hence, each of those six polyhedron is scissor equivalent to a cube if and only if $h = 1/\sqrt{2}$.

In this family $T_1(h) = P_4(h), T_2(h)$ and $T_4(h)$ are pyramids with rectangular basis. The three remaining ones : $T_3(h), T_5(h)$ and $T_6(h)$ are tetrahedron. The last one $T_6(h)$ is a Shlafli orthoscheme with lengths $(1/\sqrt{2}; 1/\sqrt{2}; h)$. As a consequence we get that a Shlafli orthoscheme with lengths $(1/\sqrt{2}; 1/\sqrt{2}; h)$ has trivial Dehn invariant if and only if $h = 1/\sqrt{2}$.

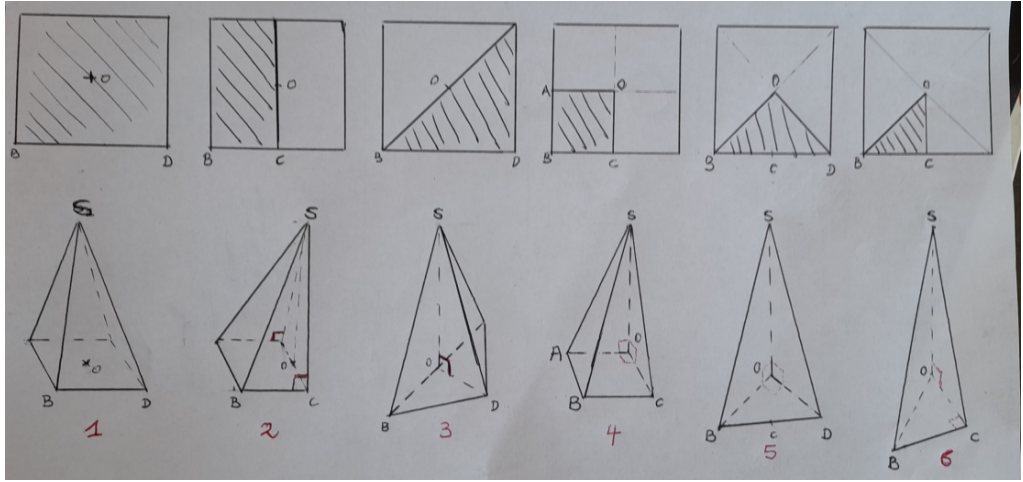

**Figure 2.**

The proof of Theorem 1, will be divide into several steps. It involves three main ingredients :

A) A theorem of Conway and Jones [3] about rational linear relations between roots of unity (see Lemma 2 below). This help to characterise rational pyramids, i.e, pyramid like $P_4(1/\sqrt{2})$ whose dihedral angles all are rational multiple of $\pi$.

B) The second ingredient is contained in Lemma 3, below. It characterises abelian Kummer extensions of number fields, when the ground field contains very few roots of the unity.

C) Finally we will conclude by using some properties of the factorisations of ideals in Dedeking domains. This will drive us to Diophantine equations with very few integral solutions, allowing to conclude.

Let's make an analogical observations about this second ingredient. In differential Galois theory, the so called "Morales-Ramis" theorem gives a strong implication : If an Hamiltonian system is "sufficiently regular" in the sens that it is "Liouville integrable" then, some associated differential Galois group must be abelian. Here, we shall meet the same idea in a geometric context. If a pyramid is "sufficiently regular" in the sens that it is "scissor equivalent to a cube", then, some associated number field extension must be abelian. This analogy is too astonishing to be not mentioned.

I solved this problem in february 2022. I then realised that this problem was solved for the pyramids $P_4(h)$ by the "anonymous user 145307" on Mathoverflow to the rubric : "Pyramid which are scissor equivalent to a cube". His work made in 2019, follows the same lines than mine modulo some important missing details that I'll present later in this development. I've tried without success to enter in contact with him. In the present redaction, I'll keep his notations, as a recognition device.

The plan of the paper will be as follow :

- In Section 2: we present the two main ingredients. But the reader which is interested to see how they apply should jump directly to Section 3.

- In Section 3 : we treat the rationality questions and prove Theorem 1, for the square basis pyramids.

- In Section 4 : we prove Theorem 1, for the pyramids with regular triangles and hexagonal basis. The two proofs could have been done simultaneously. But since they are quite involved, for sake of clarity we preferred to separate them.

- Finally, in Section 5, as a byproduct of this work, we present some non rational crystals which are scissor equivalent to cubes.

## 2 The two main ingredients of the proof

### 2.1 Linear relations between roots of the unity with rational coefficients

Here we present just a particular case of Theorem 7 of [3]. Let's consider the equation :

$$A \cos(\pi a) + B \cos(\pi b) = C,$$

where all the five variables are rational numbers. Their possible solutions are classified thanks to the following

**Lemma 2.** *Suppose that we have two rational multiple of $\pi$ angles contained into the interval $]0; \frac{\pi}{2}]$, for which some rational linear combination of their cosines is rational. Then the appropriate linear combination is proportional to*

$$A \cos\left(\frac{\pi}{2}\right) + 2 \cos\left(\frac{\pi}{3}\right) = 1.$$

*In other words up to permutation of the angles we have that $\pi a = \frac{\pi}{2}$ and $\pi b = \frac{\pi}{3}$.*

### 2.2 Abelian Kummer extensions with few roots of the unity in the ground field

Let $F$ be a field of characteristic zero and $\alpha$ a non zero element of $F$. Let's denote by $K_N(\alpha)$ the decomposition field of $X^N - \alpha$ over $F$. We shall also denote by $\mu_N$ the complex $N$-roots of the unity, and by $\text{Prim}(N)$ the subset of $\mu_N$ made of the primitive $N$-roots. We have that $K_N(\alpha) = F[u; \mu_N]$ for any root $u$ of $X^N - \alpha$. Let $G := \text{Gal}(K_N(\alpha)/F)$, be the corresponding Galois group. Let $\zeta \in \text{Prim}(N)$. Any $\sigma \in G$ act by formula of the form :

$$\sigma(\zeta) = \zeta^{a(\sigma)}, \quad \sigma(u) = u \zeta^{b(\sigma)}, \qquad (1)$$

where $(a(\sigma); b(\sigma)) \in (\mathbb{Z}/N\mathbb{Z})^* \times \mathbb{Z}/N\mathbb{Z}$. Let $\mathbb{A}_{\text{ff}}(\mathbb{Z}/N\mathbb{Z})$, be the group of affine transformations : $\Sigma(x) := a x + b$ of $\mathbb{Z}/N\mathbb{Z}$. The composition of two such formulae shows that there is an injectiv group morphism $G \to \mathbb{A}_{\text{ff}}(\mathbb{Z}/N\mathbb{Z})$. In the classical Kummer theory, that is when the ground field $F$ contains the roots of the unity, $G$ acts trivially on $\mu_N$, hence by (1), $G$ injects into the subgroups of the translations $T := \{\Sigma(x) := x + b, b \in \mathbb{Z}/N\mathbb{Z}\}$ of the affine group. Since $T \simeq \mathbb{Z}/N\mathbb{Z}$, $G$ is always abelian in this situation. Here we are interested in the context where the Kummer extension $K_N(\alpha)/F$ still is abelian but when the ground field contains very few root of unity, in the sens that we will assume that

$$F \cap \mu_N \subset \mu_2 = \{\pm 1\}.$$

Things are going to be more complicated when $N$ is an even integer. In this case we shall need two more technical "open" assumptions :

- $\mathcal{A}_1$: $\alpha$ is not a square of an element of $F$ and $-4\alpha$ is not a fourth power. That is $\alpha \notin F^2$ and $-4\alpha \notin F^4$.
- $\mathcal{A}_2$: The fields $F$ and $\mathbb{Q}[\mu_{2^n}]$ are linearly disjoint over $\mathbb{Q}$. That is $F \cap \mathbb{Q}[\mu_{2^n}] = \mathbb{Q}$, when $n \geqslant 2$.

We guess that the following result may be well known but for lake of references we prove it. There must be some links between this lemma and class field theory, but we do not manage to find them.

**Lemma 3.** *Let $N \geqslant 2$ be an integer $F$ be a field of characteristic zero such that $F \cap \mu_N \subset \mu_2 = \{\pm 1\}$. Let $\alpha$ be a non zero element of $F$. If the Kummer extension $K_N(\alpha)/F$ is abelian then :*

1. *When $N$ is odd, then $\alpha$ is a $N$-power of some element of $F$ (i.e. $\alpha \in F^N$).*
2. *When $N = 2N'$ for some odd $N'$, then $\alpha \in F^{N/2}$.*
3. *When $N = 2^n N'$, with $n \geqslant 2$, if moreover assumptions $\mathcal{A}_1$ and $\mathcal{A}_2$ hold true then*

$$\alpha = -\lambda^{N/2},$$

*for some $\lambda \in F$.*

*Conversely, those necessary conditions are sufficient for abelianity.*

In Point 3, when $N = 4N'$, observe that assumption $\mathcal{A}_2$ hold true since it is equivalent to the assumption : $F \cap \mu_4 \subset \mu_2 = \{\pm 1\}$. This last condition is satisfied since we already assumed in this case that $F \cap \mu_{4N'} \subset \mu_2 = \{\pm 1\}$.

In the present paper, the third point will be only used for $N = 4$. The proof of the remaining properties is therefore useless in the present context. Nevertheless, we develop them for later use, completeness and proper interest to Galois theory.

For its proof, we will need the following

**Theorem 4.** *Let $F$ be a field of characteristic zero and $p \geqslant 2$ be a prime number. Let $\alpha \in F$, with $\alpha \notin F^p$ (i.e. $\alpha$ is not a $p$-power of an element of $F$). Then :*

1. *$X^p - \alpha$ is irreducible over $F$.*
2. *When $p$ is odd this further implies that $X^{p^n} - \alpha$ is irreducible over $F$ for any $n \geqslant 2$.*
3. *For $p = 2$, Let $u$ belong to some extension of $F$ with $u^2 = \alpha$. Then $u$ is a square in $F(u)$ if and only if and only if $-4\alpha \in F^4$.*
4. *For any $n \geqslant 2$, $X^{2^n} - \alpha$ is irreducible over $F$ if and only if $\alpha \notin F^2$ and $-4\alpha \notin F^4$, that is if and only if $X^4 - \alpha$ is irreducible. That is if and only if assumption $\mathcal{A}_1$ is satisfied.*

This statement can be found in ([6], Theorem 9.1 p.297) and more precisely in ([7], Lemme 4.6 p.101, Theorem 4.7 p.102). Point 4, is due to Sophie Germain. These results are showing that things are easier when $N$ is odd. In light of Theorem 4, we see that in order to get the implication : $G$ abelian implies that $\alpha$ is a quasi $N/2$ power we have to add the assumption that $X^4 - \alpha$ is irreducible over the ground field.

**A general argument independent of the parity of $N$**

Let $u \in K_N(\alpha)$ be such that $u^N = \alpha$. And assume up to the end, that $G$ is abelian. Then the sub-extension $F(u)/F$ is normal, hence for any $\sigma \in G$, $\sigma(u) \in F(u)$, hence $\sigma(u)/u$ which is according to (1) a root of the unity must belong to $F[u]$.

**Descent with the odd divisors of $N$**

Let $p$ be an odd prime divisor of $N$, let $x := u^{N/p}$, then $x^p = \alpha$, hence $K_N(\alpha)$ contains $K_p(\alpha) = F[x; \mu_p]$.

If $\alpha \notin F^p$, then there exist some $\sigma \in G$ such $\sigma(x)/x = \zeta$ is a $p$ root of 1 distinct of 1. But any such $\zeta$ is primitive, hence $\mu_p \subset F(x)$, according to the previous argument. This further implies that

$$F \subset F(\mu_p) \subset F(x).$$

But the degree $d$ of $F(\mu_p)/F$ is bounded by $p - 1$, since $F(\mu_p)$ is the decomposition field of the cyclotomic polynomial : $\Phi_p(X) = X^{p-1} + \cdots + X + 1$. This degree $d$ must also divide $p = \deg(F(x)/F)$, since $X^p - \alpha$ is irreducible. Hence $d = 1$ and $F(\mu_p) = F$ which is contradictory since we assumed that $F \cap \mu_N \subset \mu_2 = \{\pm 1\}$.

As a consequence we have seen that for any odd prime divisor $p$ of $N$, $\alpha = \beta^p$ for some $\beta \in F$. Since $K_N(\alpha)$ contains $K_{N/p}(\beta)$, repeating the process we find that the abelianity of $G$ implies that

$$\alpha \in F^{N'},$$

where $N'$ is the odd part of $N$ that is $N = 2^n \times N'$ for some $n \geqslant 0$. This proves the first point of the lemma.

**Proof of Lemma 3, when $N$ is a power of two**

For $N = 2 = 2^1$, there is nothing to do, since any quadratic extension is abelian. Let's assume that $N = 2^n \geqslant 4$. In this situation, we have to prove that if assumptions $\mathcal{A}_1$ and $\mathcal{A}_2$ are satisfied, then : $K_{2^n}(\alpha)/F$ abelian implies that $\alpha = -\lambda^{2^{n-1}}$, for some $\lambda \in F$.

*Step 1* : Since $\mathcal{A}_1$ is satisfied, we know that $X^{2^n} - \alpha$ is irreducible, therefore $G$ has got a transitive action on its roots. Let $u$ be one of them, the set $\{\sigma(u)/u, \sigma \in G\}$ coincides therefore with $\mu_{2^n}$. Hence, according to the previous general argument, we get the inclusion :

$$\mu_{2^n} \subset F(u) \Rightarrow F \subset F(\mu_{2^n}) \subset F(u) = K_{2^n}(\alpha).$$

According to $\mathcal{A}_2$, we have that $\deg[F(\mu_{2^n})/F] = \deg[\mathbb{Q}(\mu_{2^n})/\mathbb{Q}] = \operatorname{card}(\mathbb{Z}/2^n\mathbb{Z})^* = 2^{n-1}$. Since $\deg[F(u)/F] = 2^n$, we get that $F(u)/F(\mu_{2^n})$ is a quadratic extension. More precisely, we claim that

$$F(\mu_{2^n}) = F(u^2). \qquad (2)$$

Indeed : Let $H$ be the subgroup of $G$ that fixes $F(\mu_{2^n})$, it is an order two subgroup of $G$. But according to (1) $H \subset T \simeq (\mathbb{Z}/2^n\mathbb{Z}; +)$ the subgroup of translations when $G$ is viewed as a subgroup of the affine group. But, the only two order subgroup of $T$ is

$$H = <\tau \colon x \mapsto x + 2^{n-1}>.$$

The action of $\tau$ on $u$ is given by $\tau(u) = u \zeta^{2^{n-1}} \Rightarrow \tau(u^2) = (u \zeta^{2^{n-1}})^2 = u^2$. Hence $u^2 \in F(\mu_{2^n})$ since it is fixed by $H$. For degree reason (2) hold true.

*Step 2* : The degree $N = 4 = 2^2$ case : According to (2), we have here $u^2 \in F(\mu_4) = F[\mathrm{i}] \neq F$, since $F \cap \mu_4 = \mu_2$. So we get an expansion of the form :

$$u^2 = f_0 + f_1 \mathrm{i} \Rightarrow \alpha = u^4 = f_0^2 - f_1^2 + 2 f_0 f_1 \mathrm{i} \Leftrightarrow \begin{cases} \alpha &= f_0^2 - f_1^2 \\ 0 &= f_0 f_1 \end{cases}.$$

If $f_1 = 0$ then $\alpha = f_0^2$ would be a square of an element of $F$. Since it is not the case from assumption $\mathcal{A}_1$, we must have that $-\alpha = f_1^2 \in F^2$.

**The reader which is only interested in the proof of Theorems 1 may jump to the next paragraph, since the only case of Lemma 3 needed for the proof of Theorem 1, is when $N = 2^n N'$, for $n \leqslant 2$.**

*Step 3* :

We now proceed by induction. We assume that $K_{2^n}(\alpha)/F$ is abelian with $2^n \geqslant 8$, and assumptions $\mathcal{A}_1$ and $\mathcal{A}_2$ are satisfied. Since $K_{2^{n-1}}(\alpha) \subset K_{2^n}(\alpha)$, the sub-extension $K_{2^{n-1}}(\alpha)/F$ is also abelian. So by induction we can assume that

$$\alpha = -\lambda^{2^{n-2}},$$

for some $\lambda \in F$. Let $u \in K_{2^n}(\alpha)$, with $u^{2^n} = \alpha$. For any primitive root $\zeta \in \mu_{2^n}$, we get that :

$$(\zeta^2 \lambda)^{2^{n-2}} = \zeta^{2^{n-1}} \times \lambda^{2^{n-2}} = -\lambda^{2^{n-2}} = (u^4)^{2^{n-2}},$$

So there exist one such primitive root $\zeta \in \operatorname{Prim}(2^n)$, such that

$$u^4 = \zeta^2 \lambda.$$

Indeed, the orbit of $u^4$ under $G$ contains $2^n/4 = 2^{n-2}$ elements which is also the cardinality of the set of the element of $\zeta^2 \lambda$ when $\zeta$ span the set of primitive roots $\operatorname{Prim}(2^n)$. According to (2), $u^2 \in F[\mu_{2^n}] = F[\zeta]$. Therefore if we set :

$$Q := u^2/\zeta \Leftrightarrow u^2 = \zeta Q,$$

This element belongs to $F[\mu_{2^n}] = F[\zeta]$ and satisfies the quadratic equation :

$$u^4 = \zeta^2 \lambda = \zeta^2 Q^2 \Rightarrow Q^2 = \lambda.$$

This is the place where assumption $\mathcal{A}_2$ is going to be again necessary. Indeed, we shall use the following :

**Proposition 5.** *Let $n \geqslant 3$ and assume that the fields $F$ and $\mathbb{Q}[\mu_{2^n}]$ are linearly disjoint over $\mathbb{Q}$. That is $F \cap \mathbb{Q}[\mu_{2^n}] = \mathbb{Q}$. Then the quadratic extensions of $\mathfrak{K}$ of $F$ contained in $F[\mu_{2^n}]$ are the three following ones : $\mathfrak{K} = F[\sqrt{D}]$ for $D \in \{-1; \pm 2\}$.*

Let us assume this property to be true. It implies that $Q$ can be expanded into the form

$$Q = x + y\sqrt{D}, \quad \text{with} \quad x, y \in F, D \in \{-1; \pm 2\}.$$

This implies that

$$Q^2 = x^2 + y^2 D + 2xy\sqrt{D} = \lambda \in F, \Rightarrow \lambda = x^2 \quad \text{or} \quad \lambda = y^2 D.$$

As a consequence, $\lambda = \pm t^2$ or $\lambda = \pm 2 t^2$ for some $t \in F$. We therefore get two possibilities :

$$\begin{cases} -\alpha = \lambda^{2^{n-2}} = (\pm t^2)^{2^{n-2}} = t^{2^{n-1}} \\ -\alpha = \lambda^{2^{n-2}} = (\pm 2 t^2)^{2^{n-2}} = 2^{2^{n-2}} \times t^{2^{n-1}} \end{cases}$$

for some $t \in F$. It remains to exclude the second possibility by showing that it implies that $K_{2^n}(\alpha)/F$ is not abelian. This will conclude the induction.

Since the primitive roots $\omega \in \text{Prim}(2^{n+1})$ are describing the $2^n$ solutions of $X^{2^n} = -1$, The roots of $X^{2^n} - \alpha = X^{2^n} + 2^{2^{n-2}} \times t^{2^{n-1}}$ are of the form:

$$u = \sqrt[4]{2} \sqrt{t} \, \omega, \quad \omega \in \text{Prim}(2^{n+1}).$$

If $K_{2^n}(\alpha)/F$ is still abelian, then the composition field : $K_{2^n}(\alpha) \cdot F[\mu_{2^{n+1}}] = K_{2^n}(\alpha)[\mu_{2^{n+1}}]$, is again abelian over $F$. It will contain $\sqrt[4]{2} \sqrt{t}$, hence all the roots of

$$X^4 - \alpha' := X^4 - 2t^2.$$

But, $\alpha' = 2t^2 \notin F^2$ otherwise $\sqrt{2}$ would be in $F$. Similarly, $-\alpha' = -2t^2 \notin F^2$ and $-4\alpha' = -2(2t)^2 \notin F^2$ otherwise $i\sqrt{2}$ would be in $F$. Therefore according to what have been proved for $N = 4$, in Step 2, the Galois group of $X^4 - \alpha' = X^4 - 2t^2$ over $F$ is not abelian. This proves what we wanted by concluding the induction procedure.

**Proof of Lemma 3, when $N$ is an even integer**

Let's decompose $N = 2^n \times N'$ for some odd integer $N'$ and $n \geqslant 1$. If $K_N(\alpha)/F$ is abelian, it contains $K_{N'}(\alpha)/F$ so there exist some $\beta$ in $F$ with

$$\alpha = \beta^{N'}.$$

But $K_N(\alpha)$ also contains the extension $K_{2^n}(\beta)/F$. Therefore $\alpha$ can be expressed into the form

$$\begin{cases} N = 2N' & \Rightarrow \alpha = \beta^{N'} = \beta^{N/2} \\ n \geqslant 2, \ N = 2^n N' \ \beta = -\lambda^{2^{n-1}} & \Rightarrow \alpha = (-\lambda^{2^{n-1}})^{N'} = -\lambda^{N/2} \end{cases}$$

Let's observe that these necessary conditions for the abelianity are also sufficient. Indeed in the three cases we get that :

1. When $N$ is odd, and $\alpha = \beta^N$, then $K_N(\alpha) = F[\mu_N]$.

2. When $N = 2N'$ and $\alpha = \beta^{N/2}$, then $K_N(\alpha) = F[\mu_N, \sqrt{\beta}] = F[\mu_N] \cdot F[\sqrt{\beta}]$ is abelian as composition of abelian extensions of $F$.

3. When $N = 2^n N'$, with $n \geqslant 2$ and $\alpha = -\lambda^{N/2}$. We again get that for any $\omega \in \text{Prim}(2^{n+1})$,

$$\alpha = (-1)^{N'}(\sqrt{\lambda})^N = (\omega^{2^n})^{N'}(\sqrt{\lambda})^N = (\omega\sqrt{\lambda})^N.$$

As a consequence $K_N(\alpha)$ is contained in the composition field $F[\mu_{2N}, \sqrt{\lambda}] = F[\mu_{2N}] \cdot F[\sqrt{\lambda}]$. It is therefore abelian.

To conclude, it remains to make the

**Proof of Proposition 5**

According to ([6], Th 1.12 p.266), that $\mathbb{Q}[\mu_{2^n}]$ and $F$ are linearly disjoint implies that

$$\mathrm{Gal}(F[\mu_{2^n}]/F) \simeq \mathrm{Gal}(\mathbb{Q}[\mu_{2^n}]/\mathbb{Q}) \simeq (\mathbb{Z}/2^n\mathbb{Z})^*.$$

In particular for $n=3$, the group $(\mathbb{Z}/8\mathbb{Z})^* = \{\bar{1}; \bar{3}; \bar{5}; \bar{7}\}$ is isomorphic to the Klein group $V_4$. Any non trivial element of this group is of order two. Hence $(\mathbb{Z}/8\mathbb{Z})^*$ contains three subgroup of order two. Therefore by the Galois correspondence the three quadratic sub field of $\mathbb{Q}[\mu_8]$ are the $\mathbb{Q}[\sqrt{D}]$ for $D \in \{-1; \pm 2\}$.

Now let's assumed that $2^n \geqslant 16$. Let $\mathfrak{K}$ be a quadratic subfield of $F[\mu_{2^n}]$. By the Galois correspondence $\mathfrak{K} = F[\mu_{2^n}]^H$ for some index two subgroup $H$ of $G = (\mathbb{Z}/2^n\mathbb{Z})^*$. Since $G/H$ is of order two, for any $\sigma \in G$, $\sigma^2 \in H$. Therefore, $H$ contains the square subgroup :

$$S := \{\sigma^2, \sigma \in G\}.$$

But as it well known there is an isomorphism $G = (\mathbb{Z}/2^n\mathbb{Z})^* \simeq (\mathbb{Z}/2\mathbb{Z}, +) \times (\mathbb{Z}/2^{n-2}\mathbb{Z}, +)$, meaning that any invertible class $\bar{k}$ mod $2^n$ can be written

$$\bar{k} = \pm \bar{5}^m \quad [2^n],$$

for some $m$ in $(\mathbb{Z}/2^{n-2}\mathbb{Z}, +)$, since the order of $\bar{5}$ in $G$ is precisely $2^{n-2}$. Therefore, $S := \{(\pm \bar{5}^m)^2, m \in (\mathbb{Z}/2^{n-2}\mathbb{Z}, +)\}$, is the cyclic subgroup generated by $\overline{25}$. Its order is $2^{n-3}$ and the index of $S$ in $G$ is $\mathrm{Card}(G/S) = 2^{n-1}/2^{n-3} = 4$.

Setting $\mathcal{M} := F[\mu_{2^n}]^S$, we have that $\mathcal{M}$ and $F[\mu_8]$ are both sub fields of $F[\mu_{2^n}]$ with the same degree equal to 4 over $F$. But the element $s = \overline{25}$ of $G$ acts on any unit root $z \in \mu_{2^n}$ by the formula $s(z) = z^{25}$. We therefore have

$$s(z) = z^{25} = z \Leftrightarrow z^{24} = 1 \Leftrightarrow z^8 = 1 \Leftrightarrow z \in \mu_8,$$

since the order of an element must divide $2^n$. So we have proved that $\mathcal{M} := F[\mu_{2^n}]^S = F[\mu_8]$. Finally by the Galois correspondence, we have that

$$S \subset H \subset G \quad \Rightarrow \quad F \subset F[\mu_{2^n}]^H = \mathfrak{K} \subset F[\mu_{2^n}]^S = F[\mu_8].$$

so any quadratic sub field of $F[\mu_{2^n}]$ is contained into $F[\mu_8]$ and belongs therefore to the quoted list.

## 3 Proof of Theorem 1 for the "egyptian pyramids"

In the first three subsections we study the Dehn invariant of the pyramids $P_n(h)$, next we prove Theorem 1 for the "egyptian pyramids" that is when $n = 4$.

### 3.1 The Dehn invariant and first condition for his vanishing

Let $h := \overline{SO}$ to be the vertical height of $P_n(h)$. Let $\theta$ and $\phi$ be the two main dihedral angles : $\theta$ is the dihedral angle along an horizontal edge belonging to the regular polygone $\mu_n$ and $\phi$ is the dihedral angle along an edge coming from $S$. See Figure 1 for an illustration. Each of these angle belongs to the real interval $]0; \pi[$. For these two characteristics angles, $\theta$ and $\phi$, elementary geometry gives

$$\tan(\theta) = \frac{h}{\cos(\pi/n)}, \quad \cos(\phi) = -\frac{h^2 \cos(2\pi/n) + \cos^2(\pi/n)}{h^2 + \cos^2(\pi/n)} \tag{3}$$

For the Dehn invariant of the general pyramid, it belongs to the Dehn group : $\mathbb{R} \otimes_\mathbb{Z} \mathbb{R}/\pi\mathbb{Z}$ and it is given by

$$\mathrm{Dehn}(P_n(h)) = n\,\xi_n \quad \text{with} \quad \xi_n := 2\sin(\pi/n) \otimes \theta + \sqrt{1+h^2} \otimes \phi$$

$\mathrm{Dehn}(P_n(h)) = n\,\xi_n = 0$ if and only if $\xi_n = 0$. There are two cases for the vanishing of $\xi_n$ :

**Case A** : Each of the two tensor $\sin(\pi/n) \otimes 2\theta$ and $\sqrt{1+h^2} \otimes \phi$ is zero. This is the case if and only if $\theta$ and $\phi$ both are rational multiples of $\pi$. That is when the pyramid $P_n(h)$ is said to be "rational".

**Case B :** Each of the two tensor $\sin(\pi/n) \otimes 2\theta$ and $\sqrt{1+h^2} \otimes \phi$ in non zero, meaning that both angles are not $\pi \mathbb{Q}$, but the sum of the tensors is zero. This implies that the ratio of lengths :

$$v := \frac{\sin(\pi/n)}{\sqrt{1+h^2}},$$

is a rational number and satisfies $0 < v < \sin(\pi/n)$. Hence,

$$\xi_n = \sqrt{1+h^2} \otimes (2v\theta + \phi) = 0 \Leftrightarrow 2v\theta + \phi \in \pi\mathbb{Q}.$$

Observe that this last condition implies Bricard's condition. See [1] and [2]. We are going to separately study these two cases by showing that **Case A** only gives the pyramid $P_4(1/\sqrt{2})$ of the cube. The study of **Case B** is more elaborated and we will show that it never happen.

## 3.2 Case A : the rationality of $P_n(h)$ for $n \in \{4; 3; 6\}$

Seeking for rational pyramids, we assume that $\theta$ and $\phi$ are in $\pi\mathbb{Q}$. These angles are related by (3). By eliminating $h$ thanks to $h = \cos(\pi/n)\tan(\theta)$, the second relation of (3) linearises into the following :

$$[1 - \cos(2\pi/n)]\cos(2\theta) + 2\cos(\phi) + [1 + \cos(2\pi/n)] = 0 \qquad (4)$$

These relations for $n=4$, $n=3$ and $n=6$ read :

$$\begin{cases} n=4: & \cos(2\theta) + 2\cos(\phi) + 1 = 0 \\ n=3: & 3\cos(2\theta) + 4\cos(\phi) + 1 = 0 \\ n=6: & \cos(2\theta) + 4\cos(\phi) + 3 = 0 \end{cases}$$

**For $n=4$,**

According to Lemma 2, we must have $\cos(2\theta) = 0 \Rightarrow \theta = \frac{\pi}{4}$ and $\phi = \frac{2\pi}{3}$. This shows that $P_4(h)$ is rational only when $h = 1/\sqrt{2}$. Originally, we made a self contained study of the equation $\cos(2\theta) + 2\cos(\phi) + 1 = 0$, thanks to the "Lucas sequences" that we investigated in [4]. This study was very similar to the one presented in Mathoverflow. We do not present it here since Lemma 2, is shorter and will also directly cover the two remaining cases :

**For $n \in \{3; 6\}$,**

According to Lemma 2, none of these two equations has got rational solutions $\theta$ and $\phi$ in $\pi\mathbb{Q}$. This shows that $P_n(h)$ is never rational for $n=3$ or $n=6$.

**Remark 6.** Here is a first justification why we only deal with pyramids $P_n(h)$ for $n \in \{4; 3; 6\}$ in this paper. Indeed, in (4), the three coefficients are rational numbers if and only if $\cos(2\pi/n) \in \mathbb{Q}$. But the degree of $\cos(2\pi/n)$ over $\mathbb{Q}$ is $\varphi(n)/2$, when $n \geqslant 3$, and $\varphi(n) = 2 \Leftrightarrow n \in \{4; 3; 6\}$. Here $\varphi$ is the Euler indicatrix function. More elaborated versions of Lemma 2, is one of the basic ingredient of [5] to get their classification of all rational pyramids of $\mathbb{R}^3$.

## 3.3 The general abelian condition for $n \in \{4; 3; 6\}$

In **Case B**, we have the simultaneous conditions :

- The two angles $\theta$ and $\phi$ are not in $\pi\mathbb{Q}$.
- There is a rational number $0 < v < \sin(\pi/n)$ such that $2v\theta_n + \phi_n \in \pi\mathbb{Q}$.

**Expressions of the exponential $e^{i\phi}, e^{i\theta}$ and $e^{i2\theta}$ in term of $v$**

Since $v = \frac{\sin(\pi/n)}{\sqrt{1+h^2}} \Rightarrow h^2 = \frac{\sin^2(\pi/n) - v^2}{v^2}$, substitution in (3) gives the following three relations

$$\begin{cases} e^{i\phi_n} &= \dfrac{-\cos(2\pi/n) - v^2 + 2i\cos(\pi/n)\sqrt{\sin^2(\pi/n) - v^2}}{1 - v^2} \\ e^{i\theta_n} &= \dfrac{\cos(\pi/n)v + i\sqrt{\sin^2(\pi/n) - v^2}}{\sin(\pi/n)\sqrt{1-v^2}} \\ e^{i2\theta_n} &= \dfrac{(1+\cos^2(\pi/n))v^2 - \sin^2(\pi/n) + 2i\cos(\pi/n)v\sqrt{\sin^2(\pi/n) - v^2}}{\sin^2(\pi/n)(1-v^2)} := \alpha_n \end{cases} \qquad (5)$$

They specialise into the following formulae in the three cases :

$$n=4:\begin{cases} e^{i\phi} &= \dfrac{-v^2+i\sqrt{1-2v^2}}{1-v^2} \\ e^{i\theta} &= \dfrac{v+i\sqrt{1-2v^2}}{\sqrt{1-v^2}} \\ e^{i2\theta} &= \dfrac{3v^2-1+2iv\sqrt{1-2v^2}}{1-v^2} := \alpha \end{cases} \qquad (6)$$

$$n=3:\begin{cases} e^{i\phi} &= \dfrac{1-2v^2+i\sqrt{3-4v^2}}{2(1-v^2)} \\ e^{i\theta} &= \dfrac{v+i\sqrt{3-4v^2}}{\sqrt{3(1-v^2)}} \\ e^{i2\theta} &= \dfrac{5v^2-3+2iv\sqrt{3-4v^2}}{3(1-v^2)} := \alpha \end{cases} \qquad n=6:\begin{cases} e^{i\phi} &= \dfrac{-1-2v^2+i\sqrt{3-12v^2}}{2(1-v^2)} \\ e^{i\theta} &= \dfrac{\sqrt{3}v+i\sqrt{1-4v^2}}{\sqrt{1-v^2}} \\ e^{i2\theta} &= \dfrac{7v^2-1+2iv\sqrt{3-12v^2}}{1-v^2} := \alpha \end{cases} \qquad (7)$$

In those three cases, if we consider the numbers fields $E := E_n = \mathbb{Q}[i\sqrt{D_n}]$, with

$$E_4 := \mathbb{Q}\left[i\sqrt{1-2v^2}\right], \quad E_3 := \mathbb{Q}\left[i\sqrt{3-4v^2}\right], \quad E_6 := \mathbb{Q}\left[i\sqrt{3-12v^2}\right],$$

they all are imaginary quadratic field, that contains $e^{i\phi}$ and $e^{i2\theta} := \alpha$ in each case. That is in order to avoid confusion, for each $n \in \{4; 3; 6\}$, $e^{i\phi_n}$ and $e^{i2\theta_n} := \alpha_n$ are in $E_n$.

In each case, from the relation

$$e^{i2v\theta} = e^{i(2v\theta+\phi)} \cdot e^{-i\phi},$$

since $2v\theta + \phi \in \mathbb{Q}\pi$, the first factor is a root of unity and so lies in an abelian cyclotomic extension of the form $\mathbb{Q}[\mu_N]$ for some $N \geq 2$. Hence $e^{i2v\theta}$ belongs to a field composition of the form

$$E[\mu_N] = E_n[\mu_N] = \mathbb{Q}[i\sqrt{D}; \mu_N],$$

which is still an abelian extension of $\mathbb{Q}$. To be more precise, let's expand $v$ as an irreducible fraction and write :

$$v := \frac{a}{b} \Rightarrow e^{i2v\theta} = e^{i2\theta\frac{a}{b}} = u^a, \quad \text{with} \quad u := e^{i2\theta/b} \Rightarrow u^b = \alpha \qquad (8)$$

Since, $e^{i2v\theta} = u^a$ we have the inclusion $E[e^{i2v\theta}] \subset E[u]$. From a Bezout relation in $\mathbb{Z}$: $as+bm=1$, we get that

$$(e^{i2v\theta})^s \times (e^{i2\theta})^m = (e^{i2v\theta})^s \times \alpha^m = e^{i2\theta/b} = u.$$

Since $e^{i2\theta} = \alpha \in E$, this implies that $E[e^{i2v\theta}] = E[u]$. But since $u^b = \alpha$, this means that the field decomposition of $X^b - \alpha$ over $E$ is a Kummer extension with abelian Galois group. With the notations of Section 1, this means that :

$$K_b(\alpha)/E,$$

is an abelian Kummer extension of number field. In order to exploit the informations contained in Lemma 3, we first need to examine the assumptions $\mathcal{A}_1$ and the condition $E \cap \mu_b \subset \mu_2$.

## 3.4 Arithmetical properties of $E = E_4 = \mathbb{Q}\left[i\sqrt{1-2v^2}\right] = \mathbb{Q}\left[i\sqrt{b^2-2a^2}\right]$

The key observation is the following technical property. We recommand the reader to skeep its proof at first reading and see how it applies in the next subsection.

**Proposition 7.** *Let's write $v = a/b$ as an irreducible fraction.*

1. *If $E$ contains some root of the unity distinct of $\pm 1$, then $E = \mathbb{Q}[i]$ and $b$ must be an odd integer.*

2. *When $b$ is an even integer, then $\alpha$ and $-\alpha$ are not square in $E$, moreover $-4\alpha \notin E^4$.*

3. *If b is an even integer, it cannot be a multiple of 4.*

Observe that Point 1 is just the condition $E \cap \mu_b \subset \mu_2$, and Point 2 is assumption $\mathcal{A}_1$. In order to prove this, we will need the following general remark that will also serve later :

**Proposition 8.** *Let $E = \mathbb{Q}[i\sqrt{D}]$ be a quadratic imaginary field. Let $w$ be a non real complex number of the form $w := \frac{\varepsilon}{\sqrt{d}}$, with $\varepsilon$ in $E$ and $d > 0$ in $\mathbb{Q}$. Set $\alpha := w^2$, then $\alpha \in E$ and*

  a) *$\alpha$ is a square in $E$ if and only if $d$ is a square in $\mathbb{Q}$.*

  b) *Assume form now that $\alpha$ is not a square in $E$, we then have the equivalence :*

$$-\alpha \in E^2 \Leftrightarrow \sqrt{\frac{D}{d}} \in \mathbb{Q}$$

  c) *Independently of the previous assumptions if $\alpha$ is any element belonging to a field $E$. If $-\alpha \notin E^2$, then $-4\alpha \notin E^2$, hence $-4\alpha \notin E^4$.*

**Proof. of Proposition 8.** Point (a) : $\alpha$ is a square in $E$ if and only if $\frac{\varepsilon}{\sqrt{d}} \in E \Leftrightarrow \sqrt{d} \in E \cap \mathbb{R} = \mathbb{Q}$.

Point (b): From : $-\alpha = (iw)^2$, we get that

$$-\alpha \in E^2 \Leftrightarrow iw = \frac{-\varepsilon}{i\sqrt{d}} \in E \Leftrightarrow i\sqrt{d} \in E.$$

And this happen if and only if $i\sqrt{d}$ is $\mathbb{Q}$-proportional to $i\sqrt{D}$. That is when $\sqrt{\frac{D}{d}} \in \mathbb{Q}$. This proves the equivalence.

Point (c) : Let $\alpha \in E$. $-4\alpha/-\alpha = 2^2 \in E^2$ therefore the two elements $-\alpha$ and $-4\alpha$ are square in $E$ simultaneously or not. □

**Proof. of Proposition 7.** Point 1. Let $\zeta$ be a root of the unity of some order $N \geqslant 3$ that belongs to $E$. Since $E/\mathbb{Q}$ is quadratic we have that $\varphi(N) \leqslant 2$. That is $N \in \{3; 6; 4\}$.

If $N = 3$ or $6$, $\mathbb{Q}[\mu_N] = \mathbb{Q}[i\sqrt{3}]$. So we would have $E = \mathbb{Q}[i\sqrt{3}] = \mathbb{Q}[i\sqrt{D}] = \mathbb{Q}\left[i\sqrt{b^2 - 2a^2}\right]$. So there will exist some integer $c$ such that $D := b^2 - 2a^2 = 3c^2$. Reducing this mod three gives $b^2 \equiv 2a^2$ [3]. Since $a$ and $b$ are prime together, in this congruence, $a$ and $b$ must be prime to three. This would imply that 2 is square mod three, which is wrong. This eliminates the possibilities $N = 3$ or $6$.

If $N = 4$. It imposes to have $E = \mathbb{Q}[i] = \mathbb{Q}[i\sqrt{D}] = \mathbb{Q}\left[i\sqrt{b^2 - 2a^2}\right]$, hence an integral solution of an equation of the form $D = b^2 - 2a^2 = c^2$. This implies that $b \equiv c$ [2] and each of these two numbers must be odd since $a$ and $b$ are prime together.

Point 2. According to (6), if we set :

$$w := e^{i\theta} = \frac{v + i\sqrt{1 - 2v^2}}{\sqrt{1 - v^2}} = \frac{a + i\sqrt{b^2 - 2a^2}}{\sqrt{b^2 - a^2}} = \frac{a + i\sqrt{D}}{\sqrt{d}} := \frac{\varepsilon}{\sqrt{d}},$$

since, $\alpha = e^{i2\theta} = w^2$, we are in the context of Proposition 8. So $\alpha \in E^2$ if and only if there exist $c \in \mathbb{N}$, with $d := b^2 - a^2 = c^2 \Leftrightarrow b^2 = a^2 + c^2$. But an even value for $b$ would imply that by reducing mod 4 that the two other variables are even (since $-1$ is not a square in $\mathbb{Z}/4\mathbb{Z}$). This again contradicts the fact that $a$ and $b$ are prime together. So when $b$ is even $\alpha \notin E^2$. Let's assume this from now. The property : $\sqrt{D/d} \in \mathbb{Q}$ means that we have an irreducible fraction :

$$\sqrt{\frac{D}{d}} = \frac{k}{l} = \sqrt{\frac{b^2 - 2a^2}{b^2 - a^2}} \Leftrightarrow \frac{a^2}{b^2} = \frac{l^2 - k^2}{2l^2 - k^2}.$$

The right denominator must be even. This forces $k$ to be even and $l$ to be odd since $k/l$ is irreducible. Hence, $2l^2 - k^2 \equiv 2 - 0 \equiv 2$ [4], which contradict the fact that the right denominator must be divisible by 4 since it is a multiple of $b^2$. Therefore, when $b$ is even $\sqrt{D/d} \notin \mathbb{Q}$ and $-\alpha \notin E^2$ and $-4\alpha \notin E^4$, thanks to Proposition 8 again.

Point 3. If $b$ is even then, thanks to Point 1 and 2 of the present proposition, $X^4 - \alpha$ is irreducible over $E$ and its field of decomposition $K_4(\alpha)$ is not abelian according to Lemma 3. Hence $K_4(\alpha)$ cannot be contained into $K_b(\alpha)$ since $K_b(\alpha)/E$ is abelian. This shows that $b$ cannot be a multiple of four. □

### 3.5 Playing with $b$-powers for the square basis pyramid case

From Proposition 7, $b$ is either odd or just two time an odd integer. According to Lemma 3 again, we know that the abelianity of $K_b(\alpha)/E$ implies that $\alpha$ is at least a $b/2$ power of an element of $E$. Indeed, if $b$ is odd, the only possible non trivial roots of the unity contained in $E$ are $\pm i$. Hence $E \cap \mu_b \subset \{\pm 1\}$, and by the first point of the lemma it exists some $\lambda \in E$ with $\alpha = \lambda^b$.

When $b = 2B$ with $B$ odd, then the Lemma gives the existence of some $\lambda \in E$ with $\alpha = \lambda^{b/2}$. In both cases, we have that
$$\alpha^2 = \lambda^b,$$
for some $\lambda \in E$. From (6)
$$\alpha = e^{2i\theta} = \frac{3a^2 - b^2 + 2ai\sqrt{b^2 - 2a^2}}{(b^2 - a^2)}.$$
Let's us set
$$z := z_4 := 3a^2 - b^2 + 2ai\sqrt{b^2 - 2a^2} \tag{9}$$

This complex number belongs to $\mathcal{O}_E$ the ring of integral elements of $E$. Its norm is $N(z) = (b^2 - a^2)^2$. Therefore we get the following relations :
$$\alpha = e^{2i\theta} = \frac{z}{(b^2 - a^2)} \Rightarrow e^{-2i\theta} = \frac{\bar{z}}{(b^2 - a^2)} \Rightarrow e^{4i\theta} = \alpha^2 = \frac{z}{\bar{z}} = \lambda^b. \tag{10}$$

In any unique factorisation domain like $\mathbb{Z}$, a relation of fraction of the form :
$$\frac{z}{z'} = \left(\frac{r}{r'}\right)^b,$$
would imply that the two elements $z$ and $z'$ will be the product of their great common divisor (g.c.d) times some $b$-power. In any Dedeking domain like $\mathcal{O}_E$ this translates into the following classical result that we include for clarity :

**Lemma 9.** *Let $R$ be a Dedekind domain, $b \geqslant 1$ some natural integer and assume that we have a relation of the form $\frac{z}{z'} = \left(\frac{r}{r'}\right)^b$, where all four variable belong to $R$. Then the g.c.d of the two ideals $(z)$ and $(z')$ is $(z; z')$, and we have a factorisation of the form :*
$$(z) = (z; z') I^b, \quad (z') = (z; z') I'^b,$$
*where $I$ and $I'$ are two co-maximal ideals of $R$.*

**Proof.** The g.c.d of two ideals is the smallest ideal that contains both ideals, hence $g.c.d\{(z); (z')\} = (z; z')$. From $z r'^b = z' r^b$, we get the factorisations of ideals
$$(z) \times (r')^b = (z') \times (r)^b. \tag{11}$$

Since $R$ is Dedekind, any ideal $J$ factorises uniquely into power of prime ideals $\mathfrak{p}$ of $R$. We will denote by $\nu_\mathfrak{p}(J)$ the power of $\mathfrak{p}$ in its factorisation. Since $(z) \subset (z, z')$, there are factorisations $(z) = (z; z') J$ and $(z) = (z; z') J'$ for some ideal $J$ and $J'$ of $R$. Since the g.c.d is the product of the maximal powers of primes $\mathfrak{p}$ both dividing $(z)$ and $(z')$, we get that $\nu_\mathfrak{p}((z; z')) = \min\{\nu_\mathfrak{p}(z); \nu_\mathfrak{p}(z')\}$. Hence $\nu_\mathfrak{p}(J) = \nu_\mathfrak{p}(z) - \min\{\nu_\mathfrak{p}(z); \nu_\mathfrak{p}(z')\}$. Therefore, for any prime, $\nu_\mathfrak{p}(J) \times \nu_\mathfrak{p}(J') = 0$ showing that $J$ and $J'$ are co-maximal. Now replacing the two above factorisations in (11) and simplification by the g.c.d gives
$$J \times (r')^b = J' \times (r)^b$$
So $\nu_\mathfrak{p}(J) + b\nu_\mathfrak{p}(r') = \nu_\mathfrak{p}(J') + b\nu_\mathfrak{p}(r)$. So if $\nu_\mathfrak{p}(J) \neq 0$ then $\nu_\mathfrak{p}(J) = b(\nu_\mathfrak{p}(r) - \nu_\mathfrak{p}(r'))$ is a multiple of $b$. Therefore there exist some $I$ such that $J = I^b$. □

Returning to (10) we have that $(z) = (z; \bar{z}) I^b$ for some ideal of $\mathcal{O}_E$. The key observation here is that any prime $\mathfrak{p}$ of $\mathcal{O}_E$ which divides the g.c.d $(z; \bar{z})$ must divide $(2) = 2\,\mathcal{O}_E$. Indeed, let $p$ in $\mathbb{Z}$ such that $(p) = \mathbb{Z} \cap \mathfrak{p}$. Since $N(z) = (b^2 - a^2)^2$ and $z + \bar{z} = 2\,(3\,a^2 - b^2)$ are in $(p)$, we get the congruences mod $p$ :

$$b^2 \equiv a^2 \quad [p], \quad 2\,(3\,a^2 - b^2) \equiv 0 \equiv 4\,a^2 \quad [p]$$

If $p \neq 2$, it will then divide both $a$ and $b$ contradicting that their g.c.d is one. So $p = 2$.

Now, taking the norm, we get that $N(z) = N((z, \bar{z})) \times N(I)^b$, but $N((z, \bar{z}))$ is a power of two. So we get a relations of the form :

$$N(z) = (b^2 - a^2)^2 = 2^k \times n^b, \tag{12}$$

where $n \in \mathbb{N}$. By eventually extracting its even part, we may assume from now that $n$ is odd.

## 3.6 Testing the equation $N(a, b) := (b^2 - a^2)^2 = 2^k \times n^b$ and conclusion

The idea is to see that (12) has got very few potential solutions. Indeed its possible solutions are

I. $(a; b) \in \{(1; 2), (1; 3), (3; 5)\}$.

II. Next for each of these possible solution we have to test the condition that $2\,v\,\theta + \phi \in \pi\,\mathbb{Q}$.

III. We will see that this never happen, therefore **Case B** never hold true. This will conclude the proof of Theorem 1 for the egyptian pyramids.

**Testing the equation $N(a, b) = (b^2 - a^2)^2 = 2^k \times n^b$**

**Proposition 10.** *The integral solutions of $N(a, b) = (b^2 - a^2)^2 = 2^k \times n^b$ with $g.c.d(a; b) = 1$ and $1 \leqslant a < b$ are the following ones :*

- *The "regular solution": $(a, b) = (1, 2)$ for which $N(1, 2) = 9 = 3^2$.*
- *The parametric solutions : $(a, b) = (2^s - 1, 2^s + 1), s \geqslant 1$, for which $N(a, b) = 2^{2s+4} \times 1^b$.*

In the last section, we shall explain why we call the solution $(a, b) = (1, 2)$ "regular".

**Proof.** If $n = 1$, then it writes $N(a, b) = (b^2 - a^2)^2 = 2^k$, hence $b - a$ and $b + a$ must powers of two. This quickly gives relations of the form $b = 2^s + 1, a = 2^s - 1$, for some $s \geqslant 1$.

If $n \neq 1$, then since it is odd we may assume that $n \geqslant 3$. We then get the inequality

$$3^b \leqslant n^b \times 2^k = (b^2 - a^2)^2 \leqslant b^4.$$

But the relation $3^b \leqslant b^4$ implies that $b \leqslant 7$. This gives very few value to test by hand and the only positive one is $(a, b) = (1, 2)$. □

In the square basis case, the condition

$$v = \frac{a}{b} = \frac{2^s - 1}{2^s + 1} \leqslant \frac{1}{\sqrt{2}},$$

implies that $s \leqslant 2 \Rightarrow s \in \{1, 2\}$. This gives the solutions $(a; b) \in \{(1; 3), (3; 5)\}$.

**Testing the condition $2\,v\,\theta + \phi \in \pi\,\mathbb{Q}$**

Multiplying by $b$ this is equivalent to $2\,\theta\,a + b\,\phi = \frac{2 k \pi}{N}$ for some natural integers $k$ and $N$. This means that

$$\Pi(a; b) = \Pi_4(a; b) := (e^{i 2 \theta})^a \times (e^{i \phi})^b = e^{i 2 k \pi / N} \in \mu_N \cap E \subset \{\pm 1; \pm i\}$$

from Proposition 7. In each of the three case, we can explicitly compute $e^{i 2\theta}$ and $e^{i \phi}$ thanks to (6) :

$$\begin{cases} (a; b) = (1; 2) & e^{i 2\theta} = e^{i \phi} = \frac{-1 + 2 i \sqrt{2}}{3} & \Pi_4(1; 2) = \left(\frac{-1 + 2 i \sqrt{2}}{3}\right)^3 = (e^{i\phi})^3 \\ (a; b) = (1; 3) & e^{i\phi} = \frac{-1 + 3 i \sqrt{7}}{8}, e^{i 2\theta} = \frac{-3 + i \sqrt{7}}{4} & \Pi_4(1; 3) = \frac{87 + 91 i \sqrt{7}}{256} \\ (a; b) = (3; 5) & e^{i\phi} = \frac{-9 + 5 i \sqrt{7}}{16}, e^{i 2\theta} = \frac{1 + 3 i \sqrt{7}}{8} & \Pi_4(3; 5) = \frac{-3617721 + 802165 i \sqrt{7}}{4194304} \end{cases} \tag{13}$$

Clearly none of these $\Pi_4(a;b)$ belongs to $\mu_4$. This concludes that **Case B**, never hold true for egyptian pyramids and the proof of Theorem 1 in this context.

# 4 Symmetric Pyramids with equilateral or hexagonal basis

## 4.1 Arithmetic properties of $E = E_n$ for $n = 3$ or $6$

According to Section 3.3,
$$K_b(\alpha)/E,$$
is an abelian Kummer extension of number field. Again, in order to use the informations contained in Lemma 3, we first need to examine the assumptions $\mathcal{A}_1$ and the condition $E \cap \mu_b \subset \mu_2$.

Here the analogue of Proposition 7 is :

**Proposition 11.** *In the triangular basis case (i.e. $n = 3$). Let's write $v = a/b$ as an irreducible fraction.*

1. *If $E_3$ contains some root of the unity distinct of $\pm 1$, then we have the following two possibilities*
   - $E_3 = \mathbb{Q}[i\sqrt{3}]$ and $a \equiv 0[3]$ and $g.c.d(b;3) = 1$.
   - $E_3 = \mathbb{Q}[i]$ and $b$ cannot be a multiple of four.

2. *The assumption : $b$ is a multiple of four implies that $\alpha$ and $-\alpha$ are not square in $E_3$, moreover $-4\alpha \notin E^4$.*

3. *If $b$ is an even integer, it cannot be a multiple of 4.*

*In the hexagonal basis case (i.e. $n = 6$). Let's write $v = a/b$ as an irreducible fraction.*

1. *If $E_6$ contains some root of the unity distinct of $\pm 1$, then we have the following two possibilities*
   - $E_6 = \mathbb{Q}[i\sqrt{3}]$ and $g.c.d(b;3) = 1$.
   - $E_6 = \mathbb{Q}[i]$ and $b$ cannot be a multiple of four

2. *The assumption : $b$ is a multiple of four implies that $\alpha$ and $-\alpha$ are not square in $E_6$, moreover $-4\alpha \notin E_6^4$.*

3. *If $b$ is an even integer, it cannot be a multiple of 4.*

We did not a find a general argument so we will have to laboriously study all the cases. We shall need the following mod 8 reduction table :

| $x \bmod [8]$ | $\pm 1; \pm 3$ | $\pm 2$ | $4; 0$ |
|---|---|---|---|
| $x^2 \bmod [8]$ | 1 | 4 | 0 |
| $4x^2 \bmod [8]$ | 4 | 0 | 0 |
| $3x^2 \bmod [8]$ | 3 | 4 | 0 |

**Proof.** of Proposition 11 Point 1. Again if $E = E_3$ or $E = E_6$ contains some non trivial root of the unity then $E = \mathbb{Q}[i]$ or $E = \mathbb{Q}[i\sqrt{3}]$.

Point 1. In the triangular basis case : $E_3 = \mathbb{Q}[i\sqrt{3}] = \mathbb{Q}\left[i\sqrt{3b^2 - 4a^2}\right] \Leftrightarrow 3b^2 - 4a^2 = 3c^2$. This forces $a$ to be a multiple of three and $b$ to be prime to three.

$E_3 = \mathbb{Q}[i] = \mathbb{Q}\left[i\sqrt{3b^2 - 4a^2}\right] \Leftrightarrow 3b^2 - 4a^2 = c^2$. This implies that $0 \equiv c^2 + b^2 \quad [4]$, hence $b$ and $c$ are even. Let's now show that this relation is impossible for $b = 4B$. This would imply that $c = 2C$ and would simplify to $12B^2 - a^2 = C^2 \Rightarrow C^2 + a^2 \equiv 0 \quad [4]$, forcing $a$ to be even which is not possible.

Point 1. In the hexagonal case : $E_6 = \mathbb{Q}[i\sqrt{3}] = \mathbb{Q}\left[i\sqrt{3b^2 - 12a^2}\right] \Leftrightarrow 3b^2 - 12a^2 = 3c^2 \Leftrightarrow b^2 = 4a^2 + c^2$. The assumption $b = 3B$ would imply that $a^2 + c^2 \equiv 0 \quad [3]$, forcing $a$ to be a multiple of three, which is not possible. Hence $g.c.d(b;3) = 1$.

$E_6 = \mathbb{Q}[\mathrm{i}] = \mathbb{Q}\left[\mathrm{i}\sqrt{3\,b^2 - 12\,a^2}\right] \Leftrightarrow 3\,b^2 - 12\,a^2 = c^2 \Rightarrow c = 3\,C \Rightarrow b^2 + 3\,C^2 = 4\,a^2$. Reduction mod 4 shows that $b$ and $C$ must have the same parity. The assumption $b = 4B$ would imply $C = 2C'$, hence $3\,C'^2 = a^2 - 4\,B^2$, with $a$ odd. But the values of $a^2 - 4\,B^2$ mod 8, belong to $\{1-4; 1-0\} = \{5; 1\}$, whose intersection with the set of three square mod 8 : $3\,C'^2 \in \{3; 4; 0\}$, is empty. Therefore, $b$ cannot be a multiple of four in this situation.

Point 2. for the triangular case : according to (7), if we set :

$$w := e^{\mathrm{i}\theta} = \frac{v + \mathrm{i}\sqrt{3 - 4\,v^2}}{\sqrt{3\,(1 - v^2)}} = \frac{a + \mathrm{i}\sqrt{3\,b^2 - 4\,a^2}}{\sqrt{3(b^2 - a^2)}} = \frac{a + \mathrm{i}\sqrt{D}}{\sqrt{d}} := \frac{\varepsilon}{\sqrt{d}},$$

since, $\alpha = e^{\mathrm{i}2\theta} = w^2$, we are in the context of Proposition 8. So $\alpha \in E^2$ if and only if there exist $c \in \mathbb{N}$, with $d := 3(b^2 - a^2) = c^2 \Rightarrow c = 3\,C \Rightarrow b^2 = a^2 + 3\,C^2$. Now let's assume that $b = 4B$. Then $a$ must be odd, and the values of $a^2 + 3\,C^2 \quad [8] \in \{1+3; 1+4; 1+0\} = \{4; 5; 1\}$ which does not meet $b^2 = 16\,B^2 \quad [8] \in \{0\}$. Therefore, $\alpha \in E^2$ implies that $b$ cannot be a multiple of four. Let's assume from now that $b = 4B$ is a multiple of four. The property : $\sqrt{D/d} \in \mathbb{Q}$ means that we have an irreducible fraction :

$$\sqrt{\frac{D}{d}} = \frac{k}{l} = \sqrt{\frac{3\,b^2 - 4\,a^2}{3\,(b^2 - a^2)}} \Leftrightarrow \frac{a^2}{b^2} = \frac{3\,l^2 - 3\,k^2}{4\,l^2 - 3\,k^2} = \frac{a^2}{16\,B^2}.$$

The right denominator must be a multiple of 16. This forces $k = 2\,K$ to be even and $l$ to be odd since $k/l$ is irreducible. Hence, $4\,l^2 - 12\,K^2 \equiv 0 \quad [16] \Rightarrow l^2 - 3\,K^2 \equiv 0 \quad [4] \Rightarrow 1^2 + K^2 \equiv 0 \quad [4]$, which is impossible. Therefore, when $b$ is a multiple of four, $\sqrt{D/d} \notin \mathbb{Q}$ and $-\alpha \notin E^2$ and $-4\,\alpha \notin E^4$, thanks to Proposition 8 again.

Point 2. For the hexagonal case : $\alpha = (e^{\mathrm{i}\theta})^2$ is a square in $E_6$ if and only if $e^{\mathrm{i}\theta} \in E_6$. But here

$$e^{\mathrm{i}\theta} = \frac{\sqrt{3}\,v + \mathrm{i}\sqrt{1 - 4\,v^2}}{\sqrt{1 - v^2}}.$$

Since $E_6$ is stable by complex conjugation $e^{\mathrm{i}\theta} \in E_6$ if and only if its real and imaginary parts are in $E$. That is $\frac{\sqrt{3}}{\sqrt{1 - v^2}} \in E \cap \mathbb{R} = \mathbb{Q}$ and $\frac{\mathrm{i}\sqrt{1 - 4\,v^2}}{\sqrt{1 - v^2}} \in \mathrm{i}\sqrt{1 - 4\,v^2}\,\mathbb{Q}$. These two conditions are not compatible since $\sqrt{3} \notin \mathbb{Q}$. Similar argument shows that $-\alpha \notin E_6^2$ in this situation.

Point 3. If $b = 4\,B$ is a multiple of four, in both cases thanks to Point 1 and 2 of the present proposition, $E \cap \mu_4 = \mu_2$, $X^4 - \alpha$ is irreducible over $E$ and its field of decomposition $K_4(\alpha)$ is not abelian according to Lemma 3. Hence $K_4(\alpha)$ cannot be contained into $K_b(\alpha)$ since $K_b(\alpha)/E$ is abelian. This shows that $b$ cannot be a multiple of four. $\square$

## 4.2 Playing with *b*-powers in the triangular and hexagonal basis cases

From Proposition 11, $b$ is either odd or just twice an odd integer. And this last proposition implies that $E \cap \mu_b \subset \mu_2$ in both cases. So according to Lemma 3 again, we know that the abelianity of $K_b(\alpha)/E$ implies that $\alpha$ is at least a $b/2$ power of an element of $E$. In the hexagonal case, in order to shorten the casuistic we will need to separate the cases :

- When $b$ is odd it exists some $\lambda \in E$ with $\alpha = \lambda^b$. This implies that $\alpha^2 = \lambda^{2b}$.

- When $b = 2\,B$ with $B$ odd, then it exists some $\lambda \in E$ with $\alpha = \lambda^{b/2}$. This implies that $\alpha^2 = \lambda^b$.

- In both cases, $\alpha^2$ is a $b$ power but in the odd case it is a $2\,b$ power.

From (7), we can write :

$$\begin{cases} e^{2\mathrm{i}\theta_3} = \frac{z_3}{3\,(b^2 - a^2)} \\ e^{2\mathrm{i}\theta_6} = \frac{z_6}{b^2 - a^2} \end{cases} \text{with} \begin{cases} z_3 := 5\,a^2 - 3\,b^2 + 2\,\mathrm{i}\,a\sqrt{3\,b^2 - 4\,a^2} \\ z_6 := 7\,a^2 - b^2 + 2\,\mathrm{i}\,a\sqrt{3\,b^2 - 12\,a^2} \end{cases} \Rightarrow \begin{cases} N(z_3) = 9\,(b^2 - a^2)^2 \\ N(z_6) = (b^2 - a^2)^2 \end{cases} \quad (14)$$

In both cases $z = z_n$ belongs to $\mathcal{O}_{E_n}$ the ring of integral elements of $E_n$. Again as in (10), we get that there is some $\lambda \in E$ such that :

$$e^{4\mathrm{i}\theta} = \alpha^2 = \frac{z}{\bar{z}} = \lambda^b. \tag{15}$$

From Lemma 9, we have a factorisation of ideals $(z) = (z; \bar{z}) I^b$ for some ideal $I$ of $\mathcal{O}_E$.

We now have to treat separately both cases, by following the same strategy as in the square basis case.

### 4.3 End of the proof of Theorem 1, for the triangular basis pyramids

Here let's show that any prime $\mathfrak{p}$ of $\mathcal{O}_E$ which divides the g.c.d $(z; \bar{z})$ must divide $(2) = 2\,\mathcal{O}_E$ or $(3) = 3\,\mathcal{O}_E$ Indeed, let $p$ in $\mathbb{Z}$ such that $(p) = \mathbb{Z} \cap \mathfrak{p}$. Since $N(z) = 9\,(b^2 - a^2)^2$ and $z + \bar{z} = 2\,(5\,a^2 - 3\,b^2)$ are in $(p)$, we get the congruences mod $p$ :

$$9\,(b^2 - a^2)^2 \equiv 0 \quad [p], \quad 2\,(5\,a^2 - 3\,b^2) \equiv 0 \quad [p]$$

If $p \neq 3$, then $b^2 \equiv a^2 \quad [p] \Rightarrow 4\,a^2 \equiv 0 \quad [p] \Rightarrow p = 2$, since $p$ cannot divide both $a$ and $b$. More precisely, there is a prime above $3\,\mathcal{O}_E$ dividing the g.c.d $(z; \bar{z})$, if and only if $a \equiv 0 \quad [3]$. Indeed, in this case, 3 divides $(5\,a^2 - 3\,b^2)$, hence $a$. Conversely if $a = 3\,A$ is a multiple of three, then $z_3 := 5\,a^2 - 3\,b^2 + 2\,\mathrm{i}\,a\sqrt{3\,b^2 - 4\,a^2}$ clearly belongs to $3\,\mathcal{O}_E$, and $(z; \bar{z}) \subset 3\,\mathcal{O}_E$.

As a consequence, $N((z; \bar{z})) = 2^k \times 3^l$ for some positive integers $k$ and $l$, with $l \geqslant 1$ when $a \equiv 0[3]$. From the factorisation of the ideals we get a relation of the form

$$N_3(a, b) := N(z) = 9(b^2 - a^2)^2 = 2^k \times 3^l \times n^b,$$

with $n \in \mathbb{N}$.

**First case : when $a \not\equiv 0 \quad [3]$**

In this case $N_3(a, b) = 9(b^2 - a^2)^2 = 2^k \times n^b$. This forces to have $n = 3\,m$ hence $9(b^2 - a^2)^2 = 2^k \times 3^b \times m^b$. This further implies that $3^b \leqslant N(a, b) < 9\,b^4 \Rightarrow b \leqslant 10$. We get a finite number of test with the conditions $a \not\equiv 0[3], g.c.d(a, b) = 1, b \leqslant 10; a^2 \leqslant \frac{3\,b^2}{4}$. The only positive one is for $(a, b) = (1, 2)$ for which $N_3(1, 2) = 3^4$.

But in this case setting again $\Pi_3(a, b) = (e^{\mathrm{i}2\theta})^a \times (e^{\mathrm{i}\phi})^b$ in order to test the condition $2\,v\theta + \phi \in \pi\,\mathbb{Q}$, we get that from (7)

$$e^{\mathrm{i}\theta_3} = e^{\mathrm{i}\phi_3} = \frac{1 + 2\,\mathrm{i}\,\sqrt{2}}{3} \Rightarrow \Pi_3(1, 2) = (e^{\mathrm{i}\phi_3})^4, \tag{16}$$

is obviously not a root of unity since $e^{\mathrm{i}\phi}$ is not.

**The second case : when $a \equiv 0 \quad [3]$ that is $a = 3\,A$**

Here $N_3(a, b) = 9(b^2 - a^2)^2 = 2^k \times 3^l \times n^b$, with $l \geqslant 1$. By eventually extracting powers, of 2 and 3 contained in $n$, we may assume from now that $n$ is prime to 6. Since $b \not\equiv 0 \quad [3]$, from: $9\,(b^2 - a^2)^2 = 2^k \times 3^l \times n^b = 9\,(b^2 - 9\,A^2)^2$, we must have $l = 2$. Hence the equation can be written :

$$(b^2 - a^2)^2 = 2^k \times n^b.$$

Taking into account the condition $a < \frac{\sqrt{3}}{2}b$, from Proposition 10, the solutions of the above equation are $(a; b) \in \{(1; 2), (1; 3), (3; 5), (7; 9)\}$. It only remains to check $(a; b) = (3, 5)$ since $a = 3\,A$. But in this case $\Pi_3(3, 5) \in \mathbb{Q}[\mathrm{i}\,\sqrt{39}]$ is not a root of the unity. This conclude the study of the triangular basis case.

### 4.4 End of the proof of Theorem 1, for the hexagonal basis pyramids

From (14), $z = z_6 = 7\,a^2 - b^2 + 2\,\mathrm{i}\,a\sqrt{3\,b^2 - 12\,a^2}$, hence for primes $p$ below divisors of the g.c.d : $(z, \bar{z})$, we have the congruences :

$$\begin{cases} b^2 \equiv a^2 & [p] \\ 2\,(7\,a^2 - b^2) \equiv 0 & [p] \end{cases} \Leftrightarrow \begin{cases} b^2 \equiv a^2 & [p] \\ 2 \times 6\,a^2 \equiv 0 & [p] \end{cases},$$

hence $p \in \{2, 3\}$ since $a \not\equiv 0 \ \ [p]$. This gives an equation

$$N(a, b) := N(z) = (b^2 - a^2)^2 = 2^k \times 3^l \times n^b,$$

where we can assume again that $n$ is prime to 6.

**First case, when $l = 0$**

Here, $N(a,b) := N(z) = (b^2 - a^2)^2 = 2^k \times n^b$. From Proposition 10, the known solutions with $\frac{a}{b} \leqslant \frac{1}{2}$ are $(a;b) \in \{(1;2), (1;3)\}$. From (7), $(a;b) = (1;2)$ gives the flat pyramid since $e^{i\theta} = 1$ in this case. We therefore exclude this option. It remains to test $(a;b) = (1;3)$.

$$e^{i2\theta} = \frac{-1 + i\sqrt{15}}{4}, e^{i\phi} = \frac{-11 + 3i\sqrt{15}}{16} \Rightarrow \Pi_6(1, 3) = \frac{-1673 + 305 i\sqrt{15}}{2^{11}}.$$

And $\Pi_6(1, 3)$ is not a root of the unity.

**Second case, when $l \geqslant 1$**

$N(a,b) = (b^2 - a^2)^2 = 2^k \times 3^l \times n^b$, where $n$ is prime to 6. Since $N(a,b)$ is a square, $k = 2K$ and $l = 2L$ are even numbers. If $n \neq 1$, then $n \geqslant 5$. And we have the inequation :

$$9 \times 5^b \leqslant N(a, b) < b^4 \Rightarrow 9 \times 5^b < b^4.$$

But this never hold true. Hence $n = 1$. And we get that

$$b^2 - a^2 = 2^K \times 3^L = (b - a)(b + a).$$

Elementary arithmetic shows that the solutions of this equation with $a < b/2$ and $g.c.d(a, b) = 1$, can be listed in the two infinite parametric following families :

- $\mathcal{F}_1 := \{(a, b) = (3^s - 2^d, 3^s + 2^d) | (s, d) \in \mathbb{N}^2, d \geqslant 1, 2^d < 3^s < 3 \times 2^d\}$,
- $\mathcal{F}_2 := \{(a, b) = (2^d - 3^s, 3^s + 2^d) | (s, d) \in \mathbb{N}^2, d \geqslant 1, 3^s < 2^d < 3^{s+1}\}$.

In each case we get that $N(z) = N(a,b) = (b^2 - a^2)^2 = 2^{2d+4} \times 3^{2s}$. Something new is needed to eliminate each of these potential possibility. In both families, $b = 3^s + 2^d$ is odd. Therefore the abelianity condition of Lemma 3, implies that there is some $\lambda \in E$ such that

$$\alpha = \lambda^b \Rightarrow \frac{z}{\bar{z}} = \lambda^{2b}.$$

Lemma 9, gives therefore a factorisation of ideal in $\mathcal{O}_E$ of the form

$$(z) = (z; \bar{z}) I^{2b}.$$

I claim that $I = \mathcal{O}_E$ must be the trivial ideal. Otherwise it will be divisible by some non trivial prime $\mathfrak{p}$. Hence $N(\mathfrak{p})^{2b}$ will divide $N(z) = 2^{2d+4} \times 3^{2s}$. Therefore, at least one of the exponents : $2d + 4$ or $2s$ will majorate $2b$. That is, we will either have

$$b = 3^s + 2^d \leqslant d + 2, \quad \text{or} \quad b = 3^s + 2^d \leqslant s.$$

But each of these inequality is wrong. Hence $I = \mathcal{O}_E$.

As a consequence, $(z) = (z; \bar{z})$. Conjugating we get that $(z) = (z; \bar{z}) = (\bar{z})$. This means that $e^{i4\theta} = \frac{z}{\bar{z}}$ is a unit of $\mathcal{O}_E$. But since $E$ is an imaginary quadratic field, the units of $\mathcal{O}_E$ are roots of the unity. Since here in **"Case B"** we assumed that $\theta \notin \pi \mathbb{Q}$, this is not possible.

This ends the proof of Theorem 1.

**Some final remarks**

Completing Remark 6, formulae (5) show that in Case B, for the vanishing of the Dehn invariant of $P_n(h)$, the two exponentials : $e^{i\phi}$ and $e^{i2\theta}$ belong to a ground field $E_n$ of the form

$$E_n = \mathbb{Q}\left[\cos\left(\frac{2\pi}{n}\right)\right][i\sqrt{D_n}],$$

with $D_n \in \mathbb{Q}\left[\cos\left(\frac{2\pi}{n}\right)\right]$. This is again the reason why we studied only the cases when $n = 4; 3; 6$, in this paper since they correspond to the simple context where $\mathbb{Q}\left[\cos\left(\frac{2\pi}{n}\right)\right] = \mathbb{Q}$.

# 5 Some crystals which are scissor equivalent to cubes

## 5.1 The three "regular pyramids" and the platonic solids

From $2\,v = \frac{2\sin(\pi/n)}{\sqrt{1+h^2}}$, we see that for the "regular value" $v = \frac{1}{2} \Leftrightarrow (a;b) = (1;2)$, the pyramid $P_n(h)$ is more regular than the other in the sens that all its $2n$ edges have the same lengths

$$c_n := \sqrt{1+h_n^2} = 2\sin(\pi/n). \qquad (17)$$

The condition $v = \frac{1}{2} = \frac{\sin(\pi/n)}{\sqrt{1+h^2}} < \sin(\pi/n)$ imposes to have $n \in \{3;4;5\}$. We therefore get only three such "regular pyramids". From (17), the corresponding pyramids are $P_4(1)$ and $P_3(\sqrt{2})$ and $P_5\left(\frac{\sqrt{5}-1}{2}\right) = P_5\left(\frac{1}{\varphi}\right)$, where $\varphi := \frac{\sqrt{5}+1}{2}$ is the gold number. This means that $P_3(\sqrt{2})$ is the regular tetrahedron and that $P_4(1)$ is the right square pyramid whose four triangular faces are equilateral. This pyramid is half of the regular octahedron. The pyramid $P_5\left(\frac{1}{\varphi}\right)$, is just a cutting of the regular icosahedron.

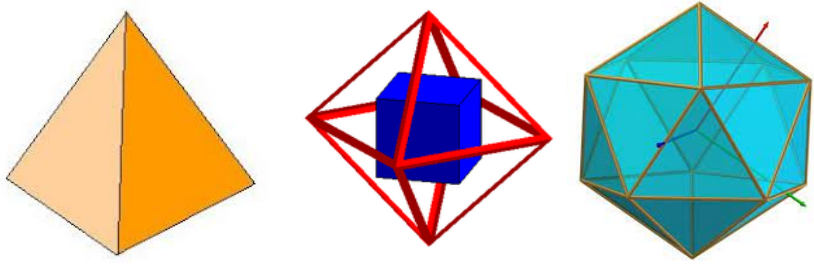

**Figure 3.**

## 5.2 Their Dehn invariants and the trivial crystals that can be deduced

According to (13) and (16), we observe the following coincidence : for $v = \frac{1}{2}$ :

$$e^{i\,2\theta_4} = e^{i\,\phi_4} = \frac{-1+2\,i\,\sqrt{2}}{3}, \quad e^{i\theta_3} = e^{i\phi_3} = \frac{1+2\,i\,\sqrt{2}}{3}.$$

Hence

$$0 < \theta_3 = \phi_3 < \pi/3, \text{ and } 2\,\pi/3 < \phi_4 = 2\,\theta_4 = \pi - \phi_3 < \pi.$$

From $\text{Dehn}(P_n(h)) = n\,\xi_n = n\,\sqrt{1+h^2} \otimes (2\,v\,\theta + \phi) = n\,\sqrt{1+h^2} \otimes (\theta + \phi)$, we get that :

$$\text{Dehn}(P_3(\sqrt{2})) = 6\,\sqrt{3} \otimes \phi_3, \quad \text{Dehn}(P_4(1)) = -6\,\sqrt{2} \otimes \phi_3.$$

This implies the two following equivalent relations :

$$\begin{cases} \sqrt{3}\,\text{Dehn}(P_4(1)) + \sqrt{2}\,\text{Dehn}(P_3(\sqrt{2})) & = 0 \\ 3\sqrt{2}\,\text{Dehn}(P_4(1)) + 2\sqrt{3}\,\text{Dehn}(P_3(\sqrt{2})) & = 0 \end{cases}. \qquad (18)$$

The second relation is just the first one multiplied by $\sqrt{6}$ but its interpretation gives something different as we shall see now.

For any polyhedron $\mathcal{P}$ and any positive real number $\lambda$, let's denote by $\lambda * \mathcal{P}$, to be the homothetic image of $\mathcal{P}$ under the ratio $\lambda$. We get that $\text{Dehn}(\lambda * \mathcal{P}) = \lambda\,\text{Dehn}(\mathcal{P})$.

**Interpretation of the relation $\sqrt{3}\,\text{Dehn}(P_4(1)) + \sqrt{2}\,\text{Dehn}(P_3(\sqrt{2})) = 0$**

Let $\mathcal{K} := \sqrt{3} * P_4(1)$ and $\mathcal{T} := \sqrt{2} * P_3(\sqrt{2})$. Each of them has edges of same size $c_\mathcal{K} = c_\mathcal{T} = \sqrt{6}$ and $\text{Dehn}(\mathcal{K}) + \text{Dehn}(\mathcal{T}) = 0$. This means that if we glue $\mathcal{K}$ and $\mathcal{T}$ together we get a polyhedron which is scissor equivalent to a cube. This is not in fact a surprise since in fact this glueing is a prism.

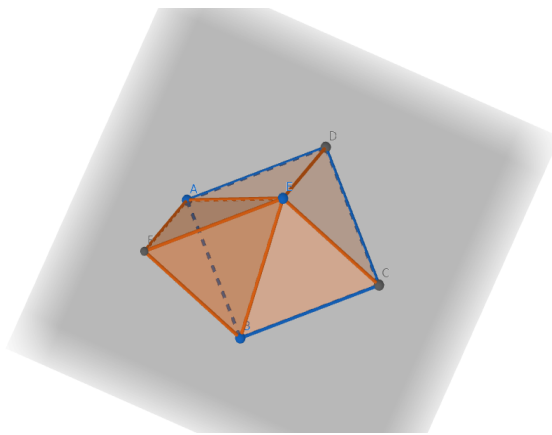

**Figure 4.**

Since $\phi_3 + \phi_4 = \pi$, we see that in the glueing edges $AE$ and $BE$ of the prism are disappearing. Something much more interesting comes now :

**Interpretation of the relation $3\sqrt{2}\,\mathrm{Dehn}(P_4(1)) + 2\sqrt{3}\,\mathrm{Dehn}(P_3(\sqrt{2})) = 0$**

Let $\mathcal{K}' := \sqrt{2} * P_4(1)$ and $\mathcal{T}' := \sqrt{3} * P_3(\sqrt{2})$. Here the sizes of each pyramid are $c_{\mathcal{K}'} = 2$ and $c_{\mathcal{T}'} = 3$ and $3\,\mathrm{Dehn}(\mathcal{K}') + 2\,\mathrm{Dehn}(\mathcal{T}') = 0$. This means that if we glue three copies of $\mathcal{K}'$ and two copies of $\mathcal{T}'$ together we get a polyhedron which is scissor equivalent to a cube. This very nice crystal can be seen in the next figure and that time it is not obvious that it is scissor equivalent to a cube.

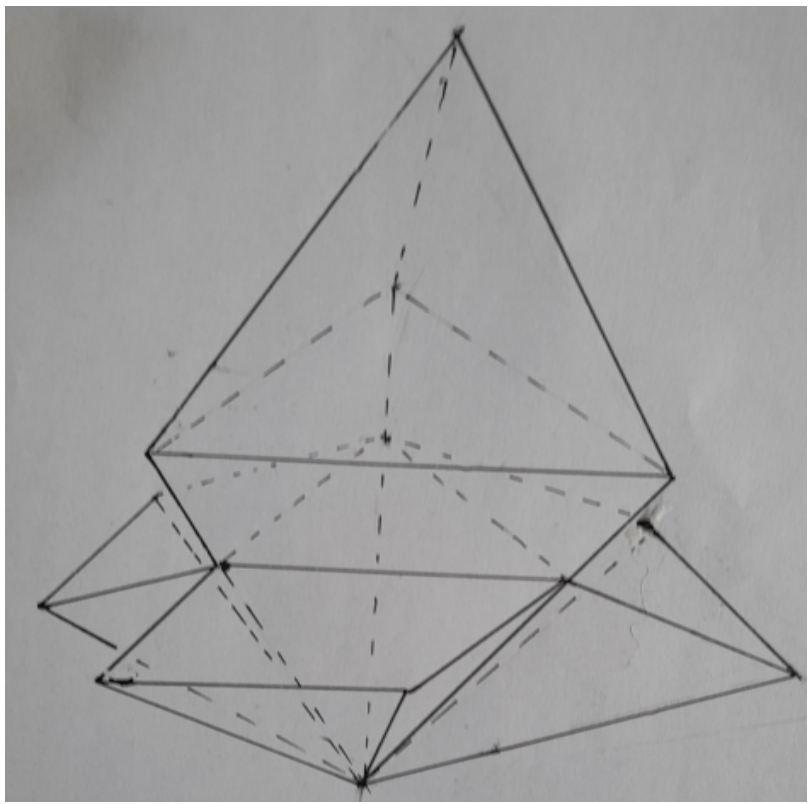

**Figure 5.**

I've used the letters $\mathcal{K}$ and $\mathcal{K}'$ for the homothetic images of $P_4(1)$ since I thought initially that they where similar to the beautiful Kheop's pyramid. This is not the case. Indeed, the Kheop's pyramid seems to be similar to $P_4\left(\frac{7\sqrt{2}}{11}\right)$, and $\frac{7\sqrt{2}}{11} \simeq 0.9$. Indeed, for the Kheop's pyramid, it seems that $\tan(\theta) = 14/11 = h\sqrt{2}$. This explains why, visually $P_4(1)$ looks very closed to the Kheop's pyramid.

For $P_5\left(\frac{1}{\varphi}\right)$, from : $\cos\left(\frac{\pi}{5}\right) = \frac{\varphi}{2} = \frac{\sqrt{5}+1}{4}$ and $\sin\left(\frac{\pi}{5}\right) = \frac{1}{2}\sqrt{\frac{5-\sqrt{5}}{2}}$ we get that

$$e^{i\phi} = \frac{-\sqrt{5}+2i}{3}, \quad e^{i\theta} = \frac{1+\sqrt{5}+(2\sqrt{5}-2)i}{\sqrt{6(5-\sqrt{5})}} \Rightarrow e^{i(\phi+\theta)} = \frac{-5\sqrt{5}-1+(4\sqrt{5}-8)i}{3\sqrt{6(5-\sqrt{5})}}$$

It can be deduced (not so easily) that $e^{i(\phi+\theta)}$ is not a root of the unity. Therefore $P_5(1/\varphi)$ is not scissor equivalent to a cube. Indeed, since $i, \sqrt{3} \in \mathbb{Q}[e^{i\pi/6}]$, $\sqrt{2}, \in \mathbb{Q}[e^{i\pi/4}]$ and $\sqrt{5}, i\sqrt{\frac{5-\sqrt{5}}{2}} \in \mathbb{Q}[e^{i\pi/5}]$ we get that $Z := e^{i(\phi+\theta)} \in \mathbb{Q}[e^{i\pi/6}, e^{i\pi/4}, e^{i\pi/5}] = \mathbb{Q}[e^{i\pi/60}]$. Therefore, if $Z$ is a root of the unity it must satisfy the relation $Z^{60} = \pm 1$. But this is not true.

## 5.3 Complexity of a polyhedron and glueing

For any non zero tensor $t \in \mathbb{R} \otimes_{\mathbb{Z}} \mathbb{R}/\pi\mathbb{Z}$, let's define its complexity $\mathrm{Comp}(t)$, to be the minimal length $l \geqslant 1$ of its decomposition as a sum of elementary tensors :

$$t = \sum_{i=1}^{l} x_i \otimes \theta_i, \quad x_i > 0, \quad \theta_i \in \mathbb{R}/\pi\mathbb{Z}.$$

For such a minimal decomposition, the family of "lengths" $\{x_1, ..., x_l\}$ is $\mathbb{Q}$ free in $\mathbb{R}$, and the familly of "angles" $\{\theta_1, ..., \theta_l\}$ is $\mathbb{Z}$ free in $\mathbb{R}/\pi\mathbb{Z}$. For any polyhedron, let's define it complexity to be the one of its Dehn invariant. Of course $\mathrm{Comp}(P) = 0$, when $P$ is scissor equivalent to a cube. When it is not, the complexity measures in some sens the "regularity" of $P$. For instance, the complexity is always bounded by the number of edges. But since $P_n(h)$ is axially symmetric we always get that $\mathrm{Comp}(P_n(h)) \leqslant 2 < 2n$. There are two cases when the complexity of $P_n(h)$ falls to one. The first one is when the ratio of edges : $2v = 2\sin(\pi/n)/\sqrt{1+h^2}$ is a rational number. The second situation is when the ratio of angles $r := \phi/\theta$ takes rational values. This can be obtain by continuously varying the height $h$. Indeed in this situation we get :

$$\mathrm{Dehn}(P_n(h)) = 2n\sin\left(\frac{\pi}{n}\right) \otimes \theta + n\sqrt{1+h^2} \otimes \phi = \left[2n\sin\left(\frac{\pi}{n}\right) + n\sqrt{1+h^2}\, r\right] \otimes \theta,$$

is of complexity one if it is not zero. I was led to this notion by meeting the strange glueing of the above subsection. Are there other such nice glueing which are scissor equivalent to cubes ? The Dehn group $\mathbb{R} \otimes_{\mathbb{Z}} \mathbb{R}/\pi\mathbb{Z}$, is naturally an infinite dimensional vector space over $\mathbb{R}$, and the above problem reduces to study the freedom of families made by the Dehn invariants of some pyramids. Is the problem simpler when we assume that all of them have complexity one ?